\documentclass[11pt]{article}
\usepackage{amsfonts,latexsym,amssymb,epsfig}
\usepackage{amsmath,amsthm,amstext}
\usepackage{fullpage,flafter}
\newcommand{\shlomo}[1]{}
\newcommand{\avner}[1]{}

\begin{document}
\bibliographystyle{plain}
\newtheorem{theorem}{Theorem}
\newtheorem{proposition}[theorem]{Proposition}
\newtheorem{corollary}[theorem]{Corollary}
\newtheorem{note}[theorem]{Note}
\newtheorem{lemma}[theorem]{Lemma}
\newtheorem{definition}[theorem]{Definition}
\newtheorem{observation}[theorem]{Observation}
\newcounter{fignum}
\newcommand{\figlabel}[1]
           {Figure \refstepcounter{fignum}\arabic{fignum}\label{#1}}
\newcommand{\ignore}[1]{}
\font\boldsets=msbm10
\def\Z{{\hbox{\boldsets  Z}}}
\def\R{{\hbox{\boldsets  R}}}
\def\half{\frac 1 2}
\def\F2{{\{0,1\}}}
\def\A{{\cal A}}
\def\K{{\cal K}} 
\def\KM{{{\cal K}^{-}}}
\def\maxk{{20000}}
\def \mixx {\mbox {mix}}
\newcommand{\rank} {\mbox {rank}}
\newcommand{\schreier} {\mbox {sc}}
\newcommand{\Oh} {\tilde{O}}
\newcommand{\G} {{\cal F}}
\newcommand{\gr}[9]{
  \fbox{$\begin{array}{lll}#1&#2&#3\\#4&#5&#6\\#7&#8&#9\end{array}$}
}
\def\hy{\hbox{-}\nobreak\hskip0pt} 
\newcommand{\SB}{\{\,}%
\newcommand{\SM}{\;{:}\;}%
\newcommand{\SE}{\,\}}%
\newcommand{\NP}{\text{\normalfont NP}}
\newcommand{\CNF}[1]{#1\text{\normalfont-CNF}}
\newcommand{\MU}{\text{\normalfont MU}}
\newcommand{\SAT}{\text{\normalfont SAT}}
\newcommand{\UNSAT}{\text{\normalfont UNSAT}}
\title{Families of unsatisfiable $k\CNF$ formulas with few occurrences per 
variable}
\author{ 
\parbox{8cm}{\centering Shlomo Hoory\\[4pt]
\small Department of Computer Science\\
\small University of Toronto\\
 Toronto, Ontario, Canada\\ 
\small\texttt{shlomoh@cs.toronto.edu}}
\parbox{8cm}{\centering Stefan Szeider\\[4pt]
\small Department of Computer Science\\
\small University of Durham\\
\small Durham, England, UK\\
\small \texttt{stefan.szeider@durham.ac.uk}}
}

\maketitle

\begin{abstract}
\sloppypar\noindent $(k,s)$-SAT is the satisfiability problem
    restricted to instances where each clause has exactly $k$
    literals and every variable occurs at most $s$ times. It is known
    that there exists a function $f$ such that for 
    $s\leq f(k)$ all $(k,s)$-SAT instances are satisfiable, but
    $(k,f(k)+1)$-SAT is already $\NP$\hy complete ($k\geq 3$).
    The best known lower and upper bounds on $f(k)$ are  $\Omega(2^k/k)$
    and $O(2^k/k^\alpha)$, where $\alpha=\log_3 4 - 1 \approx 0.26$.
    We prove that $f(k) = O(2^k \cdot \log k/k)$, which is tight up to
    a $\log k$ factor.
\end{abstract}

\section{Introduction}

We consider CNF formulas represented as sets of clauses,
where each clause is a set of literals. 
A literal is either a variable or a negated variable.
Let $k,s$ be fixed positive integers. We denote by $\CNF{(k,s)}$ the
set of formulas $F$ where  every clause
of $F$ has \emph{exactly} $k$ literals and each variable occurs in
\emph{at most} $s$ clauses of $F$. We denote the sets of satisfiable
and unsatisfiable formulas  by  $\SAT$ and $\UNSAT$, respectively.

It was observed by Tovey \cite{Tovey84} that all formulas in
$\CNF{(3,3)}$ are satisfiable, and that the satisfiability problem
restricted to $\CNF{(3,4)}$ is already $\NP$\hy complete.  This was
generalized in Kratochv\'{\i}l, et al.\ \cite{KratochST93} where it is
shown that for every $k\geq 3$ there is some integer $s=f(k)$ such that
\begin{enumerate}
\item all formulas in $\CNF{(k,s)}$ are satisfiable, and
\item $(k,s+1)$-SAT, the SAT problem
restricted to $\CNF{(k,s+1)}$, is already $\NP$\hy complete.
\end{enumerate}
The function $f$ can be defined for $k \geq 1$ by the equation
\[f(k):=\max \SB s \SM  \CNF{(k,s)}\cap \UNSAT=\emptyset \SE.\]
Exact values of $f(k)$ are only known for $k \leq 4$.
It is easy to verify that $f(1)=1$ and $f(2)=2$. 
It follows from 
\cite{Tovey84} that $f(3)=3$ and $f(k)\geq k$ in general. 
Also, by~\cite{Stribrna94}, we know that $f(4)=4$.

Upper and lower bounds for $f(k)$, $k=5,\ldots,9$, have been obtained
in~\cite{Dubois90,Stribrna94,BermanKarpinskiScott03,HoorySzeider04}.
For larger values of $k$, the best known lower bound, 
a consequence of Lov\'asz Local Lemma,
is due to Kratochv\'{\i}l et al.~\cite{KratochST93}:
\begin{eqnarray}\label{lowerbound}
f(k) \geq \left\lfloor \frac{2^k}{e k} \right\rfloor.
\end{eqnarray}
The best known upper bound, due to Savick{\'y} and Sgall~\cite{SavickySgall00},
is given by
\begin{eqnarray}\label{savupperbound}
f(k) \leq O\left(\frac{2^k}{k^\alpha}\right),
\end{eqnarray}
where $\alpha=\log_3 4 - 1 \approx 0.26$.

In this paper we asymptotically improve upon 
(\ref{savupperbound}), and show
\begin{eqnarray}\label{ourupperbound}
f(k) = O\left(\frac{2^k \log k}{k}\right).
\end{eqnarray}
Our result reduces the gap between the upper and lower bounds to a $\log k$
factor.
It turns out that the construction yielding the upper bound
(\ref{ourupperbound}) can be generalized. We present a class of $k$-CNF
formulas, that is amenable to an exhaustive search using dynamic programming.
This enables us to calculate 
upper bounds on $f(k)$ for values up to $k=\maxk$ improving upon
the bounds provided by the constructions underlying
(\ref{savupperbound}) and (\ref{ourupperbound}).

The remainder of the paper is organized as follows.
In Section~\ref{constructI:section}, we start with a simple construction
that already provides an $O(2^k \log^2 k/k)$ upper bound on $f(k)$.
In Section~\ref{constructII:section} we refine our construction
and obtain the upper bound (\ref{ourupperbound}).
In the last section we describe the more general construction and
the results obtained using computerized search.


\section {The first construction}\label{constructI:section}

We denote by $\K(x_1,\ldots,x_k)$ the complete unsatisfiable $k\CNF$ formula 
on the variables $x_1,\ldots,x_k$. This formula consists of all $2^k$ 
possible clauses. Let
$\KM(x_1,\ldots,x_k) = \K(x_1,\ldots,x_k) \setminus \{\{ x_1,\ldots,x_k \}\}$.
The only satisfying assignment for $\KM(x_1,\ldots,x_k)$ is
the all-False-assignment.
Also, for two CNF formulas $F_1$ and $F_2$ on disjoint sets of variables,
define their product $F_1 \times F_2$ as 
$\{ c_1 \cup c_2 : c_1 \in F_1 \mbox{ and } c_2 \in F_2 \}$.
Note that the satisfying assignments for $F_1 \times F_2$ are assignments
that satisfy $F_1$ or $F_2$.
In what follows, all logarithms are to the base of $2$.

\def\kl{{\lfloor k/l \rfloor}}
\def\card{\sharp}

\begin{lemma}\label{lemma1}
$f(k) 
\leq 2^k \cdot \min_{1\leq l \leq k} \left( (1-2^{-l})^{\kl} + 2^{-l}\right)$.
\end{lemma}

\begin{proof}

We prove the lemma by constructing an unsatisfiable $(k,s)\CNF$ formula $F$
where $s=2^k \cdot ( (1-2^{-l})^{\kl} + 2^{-l})$.
Let $k,l$ be two integers such that $1\leq l \leq k$, and 
let $u=\kl$ and $v=k-l \cdot u$. 
Define the formula $F$ as the union 
$F = F_0 \cup F_1 \cup \ldots \cup F_u$, where:
\begin{eqnarray*}
F_0
&=& 
\K(z_1,\ldots,z_v) \times \prod_{i=1}^u \KM(x_1^{(i)},\ldots,x_l^{(i)}),\\
F_i
&=&
\K(y_1^{(i)},\ldots,y_{k-l}^{(i)}) 
\times \{ \{ x_1^{(i)},\ldots,x_l^{(i)} \} \} 
\quad \mbox{for }i=1,\ldots,u.
\end{eqnarray*}

Therefore, $F$ is a $k\CNF$ formula with $n$ variables and $m$ clauses, where
\begin{eqnarray}
n &=& k + u \cdot (k-l)  \leq  l + k^2/l,\label{vars}\\
m &=& 2^v \cdot (2^l-1)^u + u \cdot 2^{k-l}
= 2^k \cdot \left( (1-2^{-l})^\kl + \kl \cdot 2^{-l}\right).\label{clauses}
\end{eqnarray}

To see that $F$ is unsatisfiable observe that any assignment satisfying $F_0$
must set all the variables $x_1^{(i)},\ldots,x_l^{(i)}$ to False for some $i$.
On the other hand, any satisfying assignment to $F_i$ must set at least one
of the variables $x_1^{(i)},\ldots,x_l^{(i)}$ to True.

To bound the number of occurrences of a variable note that 
the variables $z_j, y_j^{(i)}$, and $x_j^{(i)}$ occur 
$|F_0|, |F_i|$, and $|F_0|+|F_i|$ times, respectively.
Since $|F_0| = 2^v \cdot (2^l-1)^u = 2^k \cdot (1-2^{-l})^\kl$
and $|F_i|=2^{k-l}$, we get the required result.

\end{proof}

For $k \geq 4$, let $l$ be the largest integer satisfying
$2^l \leq k \cdot \log e / \log^2 k$. 
If follows that
\begin{eqnarray*}
(1-2^{-l})^{\kl} \leq \exp(-2^{-l} \cdot \kl) 
&\leq& 
\exp(-\frac{\log^2 k}{k \log e} \cdot (\frac{k}{l}+1))\\
&\leq& 
\frac{1}{\sqrt e} \cdot \exp(-\frac{\log^2 k}{l \log e})
\leq 
\frac{1}{\sqrt e} \cdot \exp(-\frac{\log k}{\log e})
=
\frac{1}{k\sqrt e},
\end{eqnarray*}
where the last inequality follows from the fact that $l \leq \log k$ for
$k \geq 4$.
Therefore, by Lemma~\ref{lemma1} there exists an unsatisfiable $k\CNF$ formula
$F$ where the number of
occurrences of variables is bounded by
\begin{eqnarray*}
2^k \cdot (\frac{1}{k\sqrt e} + \frac{2 \log^2 k}{k \log e}).
\end{eqnarray*}
It may be of interest that by (\ref{vars}) and (\ref{clauses}), 
the number of clauses in $F$ is $O(2^k \cdot \log k)$ 
and the number of variables is $O(k^2/\log k)$.

\begin{corollary}
$f(k) = O(2^k\cdot \log^2 k/k).$
\end{corollary}

\section {A better upper bound}\label{constructII:section}

To simplify the subsequent discussion, 
let us fix a value of $k$. We will only be concerned
with CNF formulas $F$ that have clauses of size at most $k$.
We call a clause of size less that $k$ an {\em incomplete} clause
and denote $F' = \{ c \in F : |c| < k \}$. 
A clause of size $k$ is a {\em complete} clause, 
and we denote $F'' = \{ c \in F : |c| = k \}$.

\begin{lemma}\label{lemma2}
$f(k) \leq 
\min\{ 2^{k-l+1} : l \in \{0,\ldots,k\} \mbox{ and } 
                   l \cdot 2^l \leq \log e \cdot (k-2l) \}$.
\end{lemma}
\begin{proof}
  
  Let $l$ be in $\{0,\ldots,k\}$, satisfying $l \cdot 2^l \leq \log e
  \cdot (k-2l)$, and set $s=2^{k-l+1}$.  We will define a sequence of
  CNF formulas, $F_0,\ldots,F_l$.  We require that (i) $F_j$ is
  unsatisfiable, (ii) $F'_j$ is a $(k-l+j)\CNF$ formula, (iii) $|F'_j|
  \leq 2^{k-l}$, and that (iv) the maximal number of occurrences of a
  variable in $F_j$ is bounded by $s$.  It follows that $F_l$ is an
  unsatisfiable $(k,s)$-CNF formula, implying the claimed upper bound.
  
  Set $k_j=k-l+j$ and $u_j = \lfloor (k-l+j)/(l-j+1)\rfloor$.  We
  proceed by induction on $j$.  For $j=0$, we define $F_0 =
  \K(x_1,\ldots,x_{k-l})$. It can be easily verified that $F_0$
  satisfies the above four requirements.  For $j>0$, assume a formula
  $F_{j-1}$ on the variables $y_1,\ldots,y_n$, satisfying the
  requirements.  We define the formula $F_j = \bigcup_{i=0}^{u_j}
  F_{j,i}$ as follows:
\begin{eqnarray}
F_{j,0} &=&
\K(z_1,\ldots,z_{k_j - u_j \cdot (l-j+1)}) \times
\prod_{i=1}^{u_j} \KM(x_1^{(i)},\ldots,x_{l-j+1}^{(i)}),\label{FPdef}\\
F_{j,i} &=&
F'_{j-1}(y_1^{(i)},\ldots,y_n^{(i)}) 
\times \{ \{ x_1^{(i)},\ldots,x_{l-j+1}^{(i)} \} \} 
\cup F''_{j-1}(y_1^{(i)},\ldots,y_n^{(i)}) 
\quad \mbox{for }i=1,\ldots,{u_j}.
\end{eqnarray}

It is easy to check that $F'_j$ is a $(k-l+j)\CNF$ formula.
To see that $F_j$ is unsatisfiable, observe that any assignment 
satisfying $F_{j,0}$, must set all the variables 
$x_1^{(i)},\ldots,x_{l-j+1}^{(i)}$ to False for some $i$. 
On the other hand, for any satisfying assignment to $F_{j,i}$, 
at least one of the variables $x_1^{(i)},\ldots,x_{l-j+1}^{(i)}$ must be set
to True.

Let us consider the number of occurrences of a variable in $F_j$.
Consider first the $y$-variables. 
These variables occur only in the $u_j$ duplicates of $F_{j-1}$ 
and therefore occur the same number of times as in $F_{j-1}$, 
which is bounded by $s$ by induction.
The number of occurrences of an $x$- or $z$-variable is
$|F'_{j-1}|+|F_{j,0}|$ or $|F_{j,0}|$ respectively.
By induction, $|F'_{j-1}| \leq 2^{k-l}$.
Also,
\begin{eqnarray*}
|F_j'|
&=&|F_{j,0}| 
= 2^{k_j - u_j \cdot (l-j+1)} \cdot (2^{l-j+1}-1)^{u_j} 
= 2^{k_j} \cdot (1-2^{-l+j-1})^{u_j}\\
&\leq& 2^{k-l+j} \cdot \exp(-2^{-l+j-1} \cdot u_j)
\leq 2^{k-l+j} \cdot \exp(-2^{-l+j-1} \cdot (k-2l)/l).
\end{eqnarray*}
Taking logarithms, we get
\begin{eqnarray*}
\log|F_{j,0}| 
&\leq& k-l+j - \log e \cdot 2^{-l+j-1} \cdot (k-2l)/l\\
&\leq& k-l+j - 2^{j-1} \leq k-l.
\end{eqnarray*}
Therefore, $F_l$ is an unsatisfiable $(k,s)\CNF$ formula for $s=2^{k-l+1}$,
as long as 
\begin{eqnarray}\label{leq}
l \cdot 2^l \leq \log e \cdot (k-2l).
\end{eqnarray}

\end{proof}

Let $l$ be the largest integer satisfying 
$2^l \leq \log e \cdot k / (2 \log k)$.
Then (\ref{leq}) holds for $k \geq 2$ 
and therefore we get the following:

\begin{corollary}
$f(k) \leq 2^k \cdot 8 \log_e k / k\quad\mbox{for }k \geq 2$.
\end{corollary}

\section{Even better upper bounds}\label{computer:section}

One way to derive better upper bounds on $f(k)$ is to generalize the 
construction of Section~\ref{constructII:section}.
To this end, we first define a special way to compose CNF formulas
capturing the essence of that construction.
\begin{definition}\label{compose:definition}
Let $F_1, F_2$ be  unsatisfiable CNF formulas that have clauses of size at 
most $k$ such that $F'_i$ is a $k_i\CNF$ formula for $i=1,2$.
Also, assume that $k_1 \leq k_2 < k$.
Then the formula $F_1 \circ F_2$ is defined as:
\begin{eqnarray*}
\left(
\bigcup_{c \in \KM(x_1,\ldots,x_{k-k_2})} F_{1,c}' \times c \cup F_{1,c}''
\right)
\cup 
F'_2 \times \{\{x_{1},\ldots,x_{k-k_2}\}\} \cup F''_2,
\end{eqnarray*}
where the formulas $F_{1,c}$ are copies of $F_1$ on distinct sets of 
variables. 
\end{definition}

It is not difficult to verify the following:
\begin{lemma}\label{compose:lemma}
Let $F_1, F_2$ be formulas as above, 
where the number of occurrences of a variable 
is bounded by $s \geq (2^{k-k_2}-1) \cdot |F'_1| + |F'_2|$
and let $G=F_1 \circ F_2$. 
Then $G$ is an unsatisfiable CNF formula
where each variable occurs at most $s$ times.
Furthermore, $G'$ is a $(k_1+k-k_2)\CNF$ formula,
and $|G'|=(2^{k-k_2}-1) \cdot |F'_1|$.
\end{lemma}

Given $k,s$, we ask whether we can obtain a $k\CNF$ formula using
the following derivation rules.
We start with the unsatisfiable formula $\{\emptyset\}$
as an axiom (this formula consists of one empty clause).
For a set of derivable formulas, one can apply one of the following rules:

\begin{enumerate}
\item
If  $F$ is a derived formula such that $s \geq
2\cdot|F'|$, 
then we  can derive
$F' \times \{\{x\},\{\overline{x}\}\} \cup F''$, where $x$ is a new 
variable.
\item
If $F_1, F_2$ are two derived formulas satisfying the conditions of 
Lemma~\ref{compose:lemma}, then we can derive the formula $F_1 \circ F_2$.
\end{enumerate}

\begin{note}
One can sometimes replace $F_1 \circ F_2$ in the second rule by a more 
compact formula $F_1 \circ' F_2$ that avoids duplicating $F_1$. 
Namely, the formula $F_1' \times \KM(x_1,\ldots,x_{k-k_2}) \cup F_1'' \cup 
F'_2 \times \{\{x_{1},\ldots,x_{k-k_2}\}\} \cup F''_2$. 
Although this can never reduce the number of occurrences of variables, 
this modification reduces the number of clauses and variables.
In the construction of Section~\ref{constructII:section},
we always use $\circ'$ instead of $\circ$.
\end{note}

\begin{figure}[h!]
  \centering
  \epsfig{file=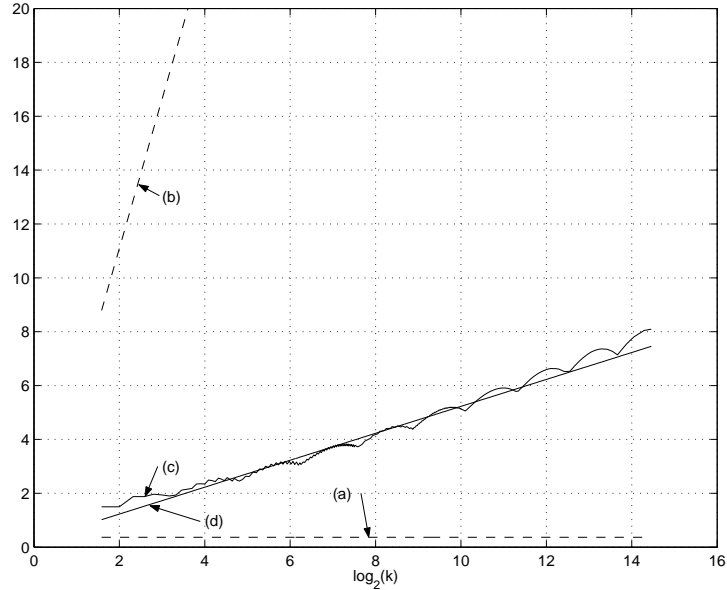,height=8cm}\\
  \caption{The bounds on $f(k) \cdot k/2^k$. 
    (a) Lower bound of Kratochv\'{\i}l et al. \cite{KratochST93},
    $1/e$.  (b) Upper bound (\ref{ourupperbound}) obtained in Section
    3 of the present paper, $8\log_e k$.  (c) Upper bound $f_2(k)
    \cdot k/2^k$, calculated by a computer program.
    (d) The line $0.5  \log(k) + 0.23$. 
  }
  \label{graph}
\end{figure}

Since any $k\CNF$ formula obtained using the above procedure is an
unsatisfiable $(k,s)\CNF$, one can define $f_2(k)$ as the maximal
value of $s$ such that no $k\CNF$ formula can be obtained using the
above procedure (clearly $f(k) \leq f_2(k)$). It turns out that the
function $f_2(k)$ is appealing from an algorithmic point of view.
Given a value for $s$, one can check if $f_2(k) \leq s$ using a simple
dynamic programming algorithm.  For all $l=0,\ldots,k-1$, the
algorithm keeps as state the minimal size of $F'$ for a derivable
unsatisfiable formula $F$ where $F'$ is an $l\CNF$ formula.  This approach
yields an algorithm that works well in practice and we were able to
calculate $f_2(k)$ for values up to $k=\maxk$ to get the results
depicted by the graph in Figure~\ref{graph}.

The computed numerical values of $f_2(k)$ seem to indicates that 
\begin{equation}
  \label{eq:conjecture}
  f_2(k) \cdot k/2^k = 0.5 \log(k) + o(\log(k))
\end{equation}
which is better than our upper bound by a constant factor of about
$11$. If (\ref{eq:conjecture}) indeed holds, then a better
analysis of the function $f_2$ may improve our upper bound by a
constant factor. However, such an approach cannot improve upon the
logarithmic gap left between the known upper and lower bounds on
$f(k)$.
 
\bibliography{f2}
\end{document}